\documentclass[12pt]{article}
\usepackage{amsfonts,amsmath,amssymb,amscd,amsthm}
\usepackage[russian]{babel}
\def\Z{{\mathbb Z}}

\long\def\comment#1\endcomment{}

\begin{document}


\newpage
\centerline{\uppercase{\bf Простое доказательство трансцендентности числа Малера}
\footnote{Это обновленная версия заметки {\it `Примеры трансцендентных чисел', Мат. Просвещение 10 (2006) 176--184.}
Заметка представлялась в 2002 г. А. Каибхановым на международной конференции Intel ISEF (США, Луисвилль),
И. Никокошевым и А. Скопенковым на Летней Конференции Турнира Городов (Россия, Белорецк) а также А. Скопенковым
в Кировской ЛМШ, Московской ОВШ и на кружках `Математический семинар', `Олимпиады и математика'.
Благодарим В. Волкова, А. Галочкина, Д. Лешко, А. Пахарева, А. Руховича и Л. Шабанова за полезные обсуждения.}
}
\smallskip
\centerline{\bf Ашум Каибханов и Аркадий Скопенков
\footnote{Поддержан грантом фонда Саймонса. Инфо: www.mccme.ru/\~{ }skopenko}
}


\bigskip
{\bf Аннотация.}
Приводятся простые доказательства трансцендентности чисел Лиувилля и Малера. 
Первое из них давно известно специалистам, а второе, по-видимому, появилось только в 2002-2003.
Доказательства доступны старшеклассникам.

\bigskip
{\bf Введение.}

Число $x$ называется {\it трансцендентным}, если оно не является корнем
уравнения
$$a_tx^t+a_{t-1}x^{t-1}+\dots+a_1x+a_0=0\qquad\text{с целыми}\qquad a_t\ne0,\ a_{t-1},\dots,a_0.$$

В университете или даже в старших классах изучается простое теоретико-множественное
доказательство существования трансцендентных чисел [KR01, Гл. 2, \S6].
Оно  не дает конкретного примера трансцендентного числа.
Приведение {\it явных} примеров трансцендентных чисел и доказательство их трансцендентности более трудно.
Первый явный пример трансцендентного числа был приведен Жозефом Лиувиллем в 1835 [KR01, Гл. 2, \S6].


\smallskip
{\bf Теорема Лиувилля.} {\it Число $\lambda=\sum\limits_{n=0}^{\infty}2^{-n!}$ трансцендентно.}


\smallskip
{\bf Общая Теорема Лиувилля.} {\it Для любых многочлена степени $t$ с рациональными коэффициентами и его иррационального корня $\alpha$ существует такое $C>0$, что для любых целых $p,q$ выполнено
$|\alpha-\dfrac pq|>Cq^{-t}$}  [KR01, Гл. 2, \S6].



\smallskip
В 1929 Курт Малер [Mah29] доказал трансцендентность следующего числа.

\smallskip
{\bf Теорема Малера.} {\it Число $\mu=\sum\limits_{n=0}^{\infty}2^{-2^n}$ трансцендентно.}

\smallskip
Эта трансцендентность не следует ни из общей теоремы Лиувилля, ни из теорем Туэ, Зигеля и Рота [KR01, Гл. 2, \S6, Fel83].
В работе [Mah29] был получен более общий результат.
Но доказательство в [Mah29] (так же как и в [Nis96]) не элементарно и длинно, ср. [Gal80].

Мы приведем короткое доказательство трансцендентности числа Лиувилля.
Это доказательство представлено в следующем пункте и известно специалистам.
К сожалению, в математических кружках обычно изучаются более сложные доказательства.

Основной результат заметки --- {\it короткое элементарное доказательство трансцендент\-ности числа Малера}
(основанное на двоичной записи).
Видимо, это доказательство не было известно до [Sko02] и [AS03, \S13.3, pp. 399-401].
Оно представлено в пункте `Доказательство теоремы Малера' и не использует остальной части заметки.
Но для удобства читателей мы представляем некоторые идеи этого доказательства в пункте `Идея простого доказательства теоремы Малера'.
Приводимые идеи дают более общий результат, см. последний пункт.

\bigskip
{\bf Доказательство теоремы Лиувилля.}

Обозначим $\lambda_s=\sum\limits_{n=0}^s2^{-n!}$.

Сначала мы докажем, что {\it число Лиувилля $\lambda$ иррационально}.
Предположим, напротив, что существует линейный многочлен
$f(x)=bx+c$ с целыми коэффициентами $b\ne0$ и $c$ такой, что $f(\lambda)=0$.
Заметим, что это уравнение имеет только один корень, значит $f(\lambda_s)\ne0$.
Мы получим противоречие из следующих неравенств для $s=|b|$:
$$2^{-s!}\ \le\ |f(\lambda_s)|\ =\ |f(\lambda)-f(\lambda_s)|\ =
\ |b|\cdot(\lambda-\lambda_s)\ <\ 2|b|\cdot2^{-(s+1)!}.$$
Первое неравенство верно, так как $f(\lambda_s)\ne0$ может быть представлена
как дробь со знаменателем $2^{s!}$.
Последнее неравенство верно, так как
$$\lambda-\lambda_s\ <\ 2^{-(s+1)!}\sum\limits_{n=0}^\infty2^{-n}\ =
\ 2\cdot2^{-(s+1)!}.$$

Сейчас мы докажем, что {\it число Лиувилля $\lambda$ не является квадратной
иррациональностью,} т.е. не является корнем квадратного уравнения
$f(x)=ax^2+bx+c=0$ с целыми коэффициентами $a\ne0$, $b$ и $c$.
Предположим, напротив, что $\lambda$ является корнем такого уравнения.
Так как квадратное уравнение имеет не более двух корней, то $f(\lambda_s)\ne0$
для достаточно больших $s$.
Тогда для достаточно больших $s$ мы получим противоречие из следующих
неравенств:
$$2^{-2s!}\ \le\ |f(\lambda_s)|\ =\ |f(\lambda)-f(\lambda_s)|\ =
\ (\lambda-\lambda_s)\cdot\left|a(\lambda+\lambda_s)+b\right| \ <
\ (2|a|\lambda+|b|)\cdot2\cdot2^{-(s+1)!}.$$
Первое неравенство верно, так как $f(\lambda_s)\ne0$ может быть представлена
как дробь со знаменателем $2^{2s!}$.
Последнее неравенство доказывается аналогично случаю линейного многочлена.

Теперь мы приведем доказательство {\it трансцендентности числа Лиувилля
$\lambda$.}
Предположим, напротив, что число $\lambda$ является корнем алгебраического
уравнения $f(x)=a_tx^t+a_{t-1}x^{t-1}+\dots+a_1x+a_0=0$ с целыми коэффициентами
$a_0, \dots, a_{t-1}, a_t\ne 0$.
Так как это уравнение имеет лишь конечное число корней, то
$f(\lambda_s)\ne0$ для достаточно больших $s$.
Тогда для достаточно больших $s$ мы получим противоречие из следующих
неравенств:
$$2^{-ts!}\ \le\ |f(\lambda_s)|\ =\ |f(\lambda)-f(\lambda_s)|\ =\
(\lambda-\lambda_s)\cdot
\left|\sum\limits_{0\le i<n\le t}a_n\lambda^{n-1-i}\lambda_s^i\right|\ <\
C\cdot2^{-(s+1)!}.$$
Первое неравенство верно, т.к. $f(\lambda_s)\ne0$ является дробью со знаменателем $2^{ts!}$.
Последнее неравенство доказывается аналогично случаю линейного многочлена.
\qed

\bigskip
{\bf Идея простого доказательства теоремы Малера.}

Продемонстрируем идею доказательства на следующем примере.
Докажем, что число
$$\nu=\sum\limits_{n=0}^{\infty}10^{-2^n}=0,11010001000000010..._{10}$$
не является {\it квадратной иррациональностью}, т.е. корнем
квадратного уравнения $x^2+bx+c=0$ с целыми коэффициентами $b$ и $c$.
(Аналогично доказывается, что $\mu$ не является квадратной иррациональностью.)
Рассмотрим десятичную запись числа $-b\nu-c$
для некоторых целых $b$ и $c$ одного знака
(случай различных знаков доказывается аналогично).
Рассмотрим ненулевые цифры в этой десятичной записи, расположенные достаточно
далеко от запятой.
Ясно, что они образуют `сгустки' около позиций с номерами $2^n$: каждый
'сгусток'
представляет число $b$.
Например для $b=-17$ мы имеем следующее:
$$17\nu-c=\dots,87170017000000170\dots017\dots.$$
А в десятичной записи числа
$$\nu^2=\sum\limits_{k,l=0}^{\infty}10^{-2^k-2^l}=
0,0121220...122020002000000012..._{10}$$
некоторые ненулевые цифры расположены около позиций с номерами $2^k+2^l$, где $k\ne l$.
Но для достаточно больших $k$ и $l$ на этих же позициях числа $-b\nu-c$ стоят нули.
Следовательно $\nu^2\ne-b\nu-c$.

\bigskip
{\bf Простое доказательство теоремы Малера.}

Пусть, напротив, $f(\mu)=a_tx^t+a_{t-1}x^{t-1}+\dots+a_1x+a_0=0$ для некоторых целых $a_t\ne0,\ a_{t-1},\dots,a_0$.
Раскрывая скобки, получим
$$\mu^q=
\left(\sum\limits_{n=0}^\infty 2^{-2^n}\right)^q=
\sum\limits_{n=0}^\infty d_n(q)2^{-n},$$
где $d_n(q)$ есть количество упорядоченных представлений числа $n$ в виде
суммы $q$ степеней двойки (не обязательно различных степеней):
$$d_n(q)=\#\{(w_1,\dots,w_q)\in\Z^q\ |\ n=2^{w_1}+\dots+2^{w_q}\text{ и }w_1,\dots,w_q>0\}.$$
Например, $d_3(2)=2$, поскольку $3=2^0+2^1=2^1+2^0$.
По определению полагаем $d_0(0)=1$.

Имеем
$$f(\mu)=\sum\limits_{n=0}^\infty d_n2^{-n},\quad\text{где}
\quad d_n:=a_td_n(t)+a_{t-1}d_n(t-1)+\dots+a_0d_n(0).$$
Ясно, что $d_n(q)=0$ тогда и только тогда, когда $n$ имеет более $q$
единиц в двоичной записи.
Для каждого $p$ положим

$\bullet$ $k=k(p):=2^{t+p}$.

$\bullet$ $m=m(p):=2^p(2^t-1)$ --- наибольшее число, меньшее $k$, для которого $d_m(t)\ne0$.

$\bullet$ $s=s(p):=2^p(2^t-1)-2^{p-1}$ --- наибольшее число, меньшее $m$, для которого $d_s(t)\ne0$.

Тогда
$$\{ 2^sf(\mu) \}\ = \ \left\{ \sum\limits_{n=0}^\infty d_n2^{s-n} \right\}\ =
\ \left\{ d_m2^{s-m}+\sum\limits_{n=k}^\infty d_n2^{s-n} \right\}.$$
Это отлично от нуля, ибо
$$\left|\sum\limits_{n=k}^\infty d_n2^{s-n}\right|\ \overset{(1)}<\ |d_m|2^{s-m}\ \overset{(2)}<\ 1/2.$$
Для доказательства неравенств (1) и (2) нам понадобится следующая лемма, которую мы докажем позже.

\smallskip
{\bf Лемма о представлении.}
{\it Количество $d_n(q)$  упорядоченных представлений числа $n$ в виде суммы
$q$ степеней двойки не превосходит $(q!)^2$.}

\smallskip
Из леммы о представлении получаем, что существует число $D=D(f)$ такое, что $|d_n|\le D$ для каждого $n$.
Следовательно неравенство (2) верно, т.к.
$|d_m|2^{s-m}\le D\cdot2^{-2^{p-1}}<1/2$ для достаточно больших $p$.
Неравенство (1) верно, поскольку для достаточно больших $p$
$$\left|\sum\limits_{n=k}^\infty d_n2^{s-n}\right|\ \le
\ D\sum\limits_{n=k}^{\infty}2^{s-n}\ =
\ D\cdot2^{s+1-k}\ =\ D\cdot2^{s-m+1-2^p} \ <\ 2^{s-m}\ \le\ |d_m|2^{s-m}.$$
Здесь последнее неравенство верно, поскольку $d_m(t)\ne0$ и $d_m(q)=0$
для $q<t$, следовательно $d_m=a_td_m(t)\ne0$.
\qed

\smallskip
{\it Доказательство леммы о представлении.} (Предложено В. Волковым.)
Индукция по $q$. Для $q=0$ имеем $d_0(0)=1\le0!^2$.
Шаг индукции вытекает из
$$d_n(q+1)\le 1+q^2d_n(q).$$
Докажем это неравенство.
Рассмотрим наборы $\vec w:=(w_1,\dots,w_{q+1})$, для которых $n=2^{w_1}+2^{w_2}+\dots+2^{w_{q+1}}$.
Наборов $\vec w$, в которых все числа различны, не более одного.
В каждом наборе $\vec w$,   в котором не все числа различны,
`объединим две равные степени двойки и поставив их на место левой из них'.
Получим некоторый набор $f(\vec w)=\vec v:=(v_1,\dots,v_q)$, для которого $n=2^{v_1}+2^{v_2}+\dots+2^{v_q}$.
Возможно, набор $f(\vec w)$ строится по набору $\vec w$ неоднозначно; выберем некоторое такое соответствие $f$.
Набор $\vec w$ получается из набора $f(\vec w)$ `разделением одной из степеней двоек на две' и `постановкой полученной
новой степени двойки в какое-то место правее исходной'.
Степень двойки для `разделения' можно выбрать $q$ способами.
`Поставить новую степень двойки в какое-то место правее исходной' можно менее, чем $q$ способами.
Поэтому у каждого набора $\vec v$, для которого $n=2^{v_1}+2^{v_2}+\dots+2^{v_q}$,
не более $q^2$ прообразов при отображении $f$.
Это доказывает нужное неравенство.
\qed

\smallskip
{\it Другое доказательство леммы о представлении.}
Для $q=0$ имеем $d_0(0)=1\le0!^2$.
Для $q\ge1$ имеем $d_n(q)=\sum\limits_{r=0}^\infty d_{n-2^r}(q-1)$.


Следовательно, при помощи индукции по $q$ получаем, что достаточно доказать следующее утверждение:

{\it для каждого $n$
существует не более $q^2$ целых $r\ge0$ таких, что $d_{n-2^r}(q-1)\ne0$.}

Докажем его.
Рассмотри двоичное представление такого числа $n$, что $d_n(q)\ne0$:
$$n=2^{w_k}+2^{w_{k-1}}+\dots+2^{w_1},\quad\text{где}\quad
w_k>w_{k-1}>\dots>w_1\ge 0\quad\text{и}\quad k\le q.$$
Обозначим $w_0=-1$.
Так как $n-2^r<0$ для $r>w_k$, то достаточно доказать, что

{\it для каждого $i=1,2,...,k$ существует не более $q$ целых $r\in [w_{i-1}+1,w_i]$ таких, что $d_{n-2^r}(q-1)\ne0$.}

Следовательно, достаточно доказать, что

{\it для каждого $i=1,2,...,k$ и $r\in[w_{i-1}+1,w_i-q]$ имеем $d_{n-2^r}(q-1)=0$.}

Чтобы доказать это утверждение, заметим, что
$2^{w_i}-2^r=2^{w_i-1}+2^{w_i-2}\dots+2^r$ есть сумма $w_i-r\ge q$
различных степеней двойки, каждая из которых больше $2^{w_{i-1}}$
и меньше $2^{w_i}<2^{w_{i+1}}$.
Поэтому число $n-2^r$ представляется в виде суммы более чем $q-1$ различных степеней двойки.
Значит, $d_{n-2^r}(q-1)\ne0$.
\qed

\bigskip
{\bf Обобщение.}

Аналогично доказывается, что $\sum\limits_{n=0}^{\infty}b_n2^{-2^n}$ трансцендентно для любой ограниченной
последовательности $\{b_n\}$ целых
чисел, среди которых бесконечно много ненулевых.

Предположим теперь, что дана строго возрастающая последовательность $\{a_n\}$ натуральных чисел.

Пусть дано натуральное число $q$.
Натуральное число $m$ называется {\it $q$-представимым}, если $m$ можно представить в виде суммы не более чем
$q$ членов последовательности $\{a_i\}$: $m=a_{i_1}+a_{i_2}+\dots+a_{i_q}$.
Эти члены не обязательно различны (например, число $m=a_2+a_2$ является 2-представимым).

Последовательность $\{a_i\}$ называется {\it $q$-разреженной}, если для любого
целого $M$ существует три последовательных $q$-представимых числа $a$, $b$ и
$c$, промежутки между которыми больше $M$ (т.е. таких, что $b-a>M$ и $c-b>M$).

Последовательность $\{a_i\}$ называется {\it разреженной}, если она $q$-разреженная для любого $q$.

Например, последовательность $a_i=i$ всех целых положительных чисел является
1-представимой, следовательно, эта последовательность не 1-разрежена и тем более не разрежена.
При доказательстве теоремы Малера было доказано, что последовательность $2^i$ является разреженной.

Последовательность $\{a_n\}$ натуральных чисел называется {\it $q$-рыхлой},
если количество способов представления числа $n$ в виде суммы $q$ членов
этой последовательности (не обязательно различных) не превосходит некоторой
константы, не зависящей от $n$ (но, возможно, зависящей от $q$).
Сформулируем это немного по-другому.
Для любых натуральных $q$ и $n$ обозначим через $d_n(q)$ количество
представлений числа $n$ в виде суммы $q$ членов последовательности $\{a_i\}$
(не обязательно различных членов).
Другими словами, $d_n(q)$ --- это число решений уравнения
$n=a_{i_1}+\dots+a_{i_q}$ с $i_1\le\dots\le i_q$.
Последовательность $\{a_i\}$ называется {\it $q$-рыхлой}, если существует число
$C_q$ такое, что $d_n(q)<C_q$ при любом $n$.

Последовательность $\{a_i\}$ называется {\it рыхлой}, если она является
$q$-рыхлой для любого натурального $q$.

Например последовательность $a_i=i$ не рыхлая и даже не 2-рыхлая.
Действительно, натуральное число $n$ имеет не менее $n/2-1$ представлений в виде суммы двух натуральных чисел, т.о. $d_n(2)\ge n/2-1$.
Из леммы о представлении вытекает, что последовательность $2^i$ является рыхлой.

Аналогично доказательству теоремы Малера можно доказать следующий результат.

\smallskip
{\bf Теорема.} {\it Если последовательность $\{a_n\}$ рыхлая и разреженная, то число
$\sum\limits_{n=0}^\infty b_n2^{-a_n}$ трансцендентно для любой ограниченной последовательности $\{b_n\}$ целых чисел.}


\smallskip
Было бы интересно выяснить, трансцендентны ли следующие числа (например, применив приведенную
теорему или метод ее доказательства):
$$\sum\limits_{n=0}^{\infty}n2^{-2^n},\quad \sum\limits_{n=0}^{\infty}2^{n-2^n},\quad
\sum\limits_{n=0}^{\infty}2^{-[1,1^n]},\quad \sum\limits_{n=0}^{\infty}2^{-f_n},$$ 
где $f_{n+2}=f_{n+1}+f_n$, $f_0=f_1=1$ --- последовательность Фибоначчи.

\smallskip
Было бы интересно обобщить наше доказательство трансцендентности числа Малера
и теоремы о редкости до достаточного условия трансцендентности числа,
включающего два
приближения этого числа.

\bigskip
\centerline{\bf Литература.}

[AS03] Jean-Paul Allouche and Jeffrey Shallit, Automatic Sequences: Theory, Applications, Generalizations, Cambridge University Press, Cambridge, 2003.

[Sko02] А. Скопенков, Примеры трансцендентных чисел,
\linebreak
http://www.turgor.ru/lktg/2002/problem5.ru/index.php.

[Fel83] Н. Фельдман, Алгебраические и трансцендентные числа, Квант, N7 (1983), 2--7.

[Gal80] А. Галочкин, О мере трансцедентности значений функций, удовлетворяющих некоторым
функциональным уравнениям, Мат. Заметки, 27:2 (1980), 175--183.

[KR01] Р. Курант и Г. Роббинс, Что такое математика? МЦНМО, Москва, 2001.

[Mah29] K. Mahler, Arithmetische Eigenschaften der L\"osungen einer Klasse von Funktional\-gleichungen,
Mathematische Annalen, 1 (1929) 342--366.

[Nis96] K. Nishioka, Mahler Functions and Transcendence, Lecture Notes in Math. {\bf 1631} (1996), Springer-Verlag, Berlin-New York.

\end{document}